\documentclass[12pt,a4paper]{amsart}

\usepackage[margin=3cm]{geometry}
\usepackage{graphicx}
\usepackage{verbatim}
\usepackage{amssymb}
\usepackage{amsbsy}
\usepackage{amscd}
\usepackage{amsmath}
\usepackage{amsthm}
\usepackage{amsxtra}
\usepackage{latexsym}
\usepackage[mathscr]{eucal}
\usepackage{xcolor}

\usepackage[T1]{fontenc}
\usepackage{tgcursor}

\usepackage{tikz}
\usetikzlibrary{arrows}

\usepackage[hidelinks]{hyperref}

\theoremstyle{plain}
\newtheorem{thm}{Theorem}[section]
\newtheorem{lemm}[thm]{Lemma}
\newtheorem{prop}[thm]{Proposition}
\newtheorem{cor}[thm]{Corollary}

\theoremstyle{definition}
\newtheorem{dfn}[thm]{Definition}
\newtheorem{exa}[thm]{Example}

\theoremstyle{remark}
\newtheorem{rmk}[thm]{Remark}

\newcommand{\fg}{\mathfrak g}

\DeclareMathOperator{\End}{End}

\newcommand{\hu}{\hat{u}}
\newcommand{\hchi}{\hat{\chi}}

\newcommand{\ZZ}{\mathbb Z}

\newcommand{\CC}{\mathbb{C}}
\newcommand{\cM}{\mathcal{M}}

\newcommand{\Sl}{\mathscr{S}}
\newcommand{\Sf}{\Sl_f}

\newcommand{\fsl}{\mathfrak{sl}}

\newcommand{\sub}{\mathrm{sub}}
\newcommand{\fin}{\mathrm{fin}}

\newcommand{\fh}{\mathfrak{h}}

\newcommand{\Jac}{\mathrm{Jac}}

\newcommand{\cW}{\mathcal W}
\newcommand{\Id}{\mathrm{Id}}
\newcommand{\cA}{\mathcal A}
\newcommand{\GP}{\Gamma_+}
\newcommand{\GPz}{\Gamma_{+,0}}
\newcommand{\cH}{\mathcal H}

\newcommand{\cV}{\mathcal V}
\newcommand{\wJ}{\widehat J}
\newcommand{\wS}{\widehat S}
\newcommand{\fzp}{\mathfrak{z}_p}
\newcommand{\AAA}{\mathbb{A}}
\title{Filtered \(R\)-Matrices: PBW Quantization and Integrable Systems}
\author{Wenda Fang}
\address{Research Institute for Mathematical Sciences, Kyoto University,
 Kyoto 606-8502 JAPAN}
\email{wenda@kurims.kyoto-u.ac.jp}
\thanks{}

\begin{document}

\maketitle

\begin{abstract}
We introduce filtered $R$-matrices for nonlinear Lie algebras and
nonlinear Lie conformal algebras and prove a splitting theorem in
each setting. For the subregular Slodowy slice in $\mathfrak{sl}_3$,
we construct finite and chiral PBW quantizations of non-splitting
$R$-deformations and obtain commuting quantum Hamiltonians for the
RDW system. At a suitable specialization, the quantum RDW Hamiltonian
yields the quantum spectral equation of the periodic $A_2$ Toda chain.
The corresponding irregular $\mathrm{SL}_3$-oper
supports the conjectural link between the RDW system and wild
Hitchin system.
\end{abstract}

\section{Introduction}

\subsection*{Motivation and the nonlinear problem}

The classical $R$-matrix formalism provides a systematic approach
to the construction of Poisson structures and integrable systems.
Its extension to Lie conformal algebras was developed by the author
in~\cite{Fang25}, leading to the generalized
Adler--Kostant--Symes (GAKS) scheme.
Motivated by the GAKS scheme, we ask when the $R$-deformation
of the bracket of a nonlinear Lie algebra again defines a
nonlinear Lie algebra. In particular, under what conditions
does a graded splitting give rise to a filtered $R$-matrix?
We study the same questions for nonlinear Lie conformal algebras. The theory of nonlinear Lie algebras (NLAs) and nonlinear
Lie conformal algebras (NLCAs), developed by De Sole
and Kac~\cite{DSK05}, provides a natural algebraic
framework for these questions.
Unlike ordinary Lie algebras, a nonlinear Lie algebra
is equipped with a bracket taking values in the
tensor algebra of its generating space.
Moreover, its Jacobi identity is formulated modulo
filtered relation spaces determined by the nonlinear bracket.

Consequently, given a nonlinear Lie algebra $(V,\mu)$
and a degree preserving linear operator $R:V\to V$,
the deformed bracket
\[
\mu_R(a,b)
=
\mu(Ra,b)+\mu(a,Rb)
\]
does not automatically define another nonlinear Lie algebra.
The essential difficulty is that the relation spaces
associated with $\mu$ and $\mu_R$ are generally different.

Thus, the PBW theorem alone does not establish the
admissibility of an $R$-deformation.
One must first verify that the deformed bracket satisfies
the nonlinear Jacobi condition with respect to its own
filtered relation spaces.

We prove a splitting theorem in both settings and construct
finite and chiral PBW quantizations for explicit $R$-deformations
associated with the subregular Slodowy slice in $\mathfrak{sl}_3$.

\newpage
\subsection*{Main results}
\paragraph{Filtered $R$-matrices and splitting constructions.}
We introduce the notion of a filtered $R$-matrix
for nonlinear Lie algebras and its chiral analogue
for nonlinear Lie conformal algebras. Our main structural results give sufficient conditions for the
classical splitting construction to extend to nonlinear Lie algebras and nonlinear Lie conformal algebras.

\begin{thm}[Splitting theorem]
Let $(V,\mu)$ be a $\GP$-graded nonlinear Lie algebra with a
graded decomposition $V=V_+\oplus V_-$. Suppose that, for some
fixed $\delta\in\GP$ and all $p,q\in\GP$,
\[
\mu(V_\pm[p],V_\pm[q])
\subseteq T(V_\pm)[p+q-\delta].
\]
Then $R=P_+-P_-$ is a filtered $R$-matrix, where $P_\pm$ are
the projections onto $V_\pm$.

The same conclusion holds for a nonlinear Lie conformal algebra
when $V_\pm$ are $\CC[\partial]$-submodules and the displayed
condition holds coefficientwise for the $\lambda$-bracket.
\end{thm}

In the proof, we use the homogeneity assumption to compare the
filtered relation spaces of the original and rescaled brackets
on each summand. The PBW theorems of De Sole and Kac then yield
filtered associative algebras in the finite case and filtered
vertex algebras in the chiral case. Their associated graded
algebras carry the corresponding Poisson and Poisson vertex
algebra structures, respectively.

For the nonlinear Lie algebra associated with the subregular
Slodowy slice in $\mathfrak{sl}_3$, we classify all
degree-preserving splitting-type $R$-matrices
(Theorem~\ref{thm:subW3SP}).

\paragraph{Nonsplitting $R$-matrices and PBW quantization.}

The preceding splitting theorems provide a general
source of filtered $R$-matrices.
However, the theory is not restricted to
splitting operators. To investigate $R$-deformations beyond the splitting
construction, we consider finite $\mathcal W$-algebras
and the Poisson geometry of Slodowy slices. For the subregular nilpotent element of
$\mathfrak{sl}_3$, we construct families of
degree preserving filtered $R$-matrices,
including nonsplitting examples.
We establish the corresponding nonlinear Lie algebra
structures and obtain PBW quantizations of the
$R$-deformed Poisson algebras.

We further investigate the chiralization of these Poisson structures.
For the first-order differential $R$-deformations $R_4$ and $R_7$
defined in Section~\ref{sec:PBWPVA}, we construct explicit PBW
quantizations on the arc space of the subregular Slodowy slice,
retaining $\chi$ as an independent central differential generator.
For the $d=a$ branch of the mixed family described in
Theorem~\ref{thm:PBWQRDW}, we also construct an explicit chiral
PBW quantization.
A notable feature of this construction is the
appearance of quantum correction terms depending
nonlinearly on the parameters of the $R$-matrix.
The resulting quantization is valid throughout this branch, including its degenerate
specializations.

\paragraph{Classical and quantum integrability.}

The classical $R$-matrix formalism is closely related
to the construction of commuting Hamiltonians.
Let $(\fg,[\cdot,\cdot],(\cdot,\cdot))$ be a quadratic
Lie algebra and let $R:\fg\to\fg$ be a classical
$R$-matrix. Denote by $\fzp(\fg)\subset S(\fg)$ the Poisson
center with respect to the standard Lie--Poisson
bracket $\{\cdot,\cdot\}$.
The classical AKS construction implies that
\[
\{F,G\}_R=0
\qquad
\text{for all }F,G\in\fzp(\fg),
\]
thereby providing Poisson commuting Hamiltonians
for the $R$-deformed Lie--Poisson structure. The generalized AKS scheme extends this approach
to Lie conformal algebras and Poisson vertex
algebras~\cite{Fang25}. Its applications to nonlinear Poisson structures, in particular
those arising from classical $\mathcal W$-algebras, were further
investigated in~\cite{Fang26}. One example introduced there is the
classical RDW system. At $\chi=0$, its generic spectral curve has
genus two and admits an order-three automorphism. The induced action
on its Jacobian relates the RDW system to the generalized Kummer
surfaces studied by Bertin and Vanhaecke~\cite{BV94}. We brief review RDW system and discuss its property.

A PBW quantization of a Poisson algebra leaves open whether its
commuting Hamiltonians admit commuting quantum lifts. For the RDW
system introduced in~\cite{Fang26}, we construct such lifts of $S$
and $J$ on the mixed branch $d=a$. In the chiral quantization, we
obtain a commuting pair of quantum local Hamiltonians. The normally
ordered squares of their densities, however, need not define
commuting local Hamiltonians.

We next consider the spectral problem for the quantum RDW
Hamiltonian at a distinguished specialization of the mixed branch. 
\begin{thm}[Quantum RDW--Toda correspondence]
Consider the quantum RDW system with $a=b=0$ and $p=r=1$,
and fix the central characters $\hat\chi=c$ and $\hat J=j$,
where $c\in\CC$ and $j\in\CC^\times$.
The reduced Hamiltonian admits a differential-operator realization
in which its eigenvalue equation takes the form
\[
\left[
D_z^3+\left(3c+\frac{\epsilon^2}{8}\right)D_z
+z+\frac{9j^2}{4z}-6e
\right]\psi=0,
\qquad D_z=\epsilon z\partial_z,\quad \epsilon=-4i\hbar,
\]
where $e$ is the spectral parameter.
After rescaling $z$ and identifying the parameters, this equation
coincides with the quantum spectral equation of the periodic
$A_2$ Toda chain in~\cite[Eq.~(3.2)]{BRT25}.
It defines an $\mathrm{SL}_3$-oper on $\mathbb G_m$ over
$\CC((\epsilon))$ with irregular singularities at $0$ and $\infty$.
\end{thm}

This oper supports the conjectural link between the RDW system
and wild Hitchin geometry.

On the classical side, for $c=0$ and fixed $j\ne0$, the RDW spectral
curves form the $u=0$ subfamily of Seiberg--Witten curves for pure
$\mathcal N=2$ $\mathrm{SU}(3)$ gauge theory~\cite[Eq.~(4.4)]{KLT95}.

\subsection*{Relation to previous work and outlook}

The present paper continues the generalized AKS
program initiated in~\cite{Fang25} and extended
to nonlinear Poisson structures and RDW integrable
systems in~\cite{Fang26}. For the mixed branch $d=a$, our construction also gives a comparison
with the critical level subregular $\mathcal W$-algebra of
$\mathfrak{sl}_3$. When $ap\Delta\ne0$, the vertex algebra before
the $\hbar$-rescaling is isomorphic to
$\mathcal W^{-3}(\mathfrak{sl}_3,f_{\mathrm{sub}})$
(see Remark~\ref{rmk:mixPBA}).

Ambrosio, Carnovale, Esposito, and Topley
\cite[Theorem~1.2]{ACET24}
study universal filtered quantizations of Poisson deformations
of nilpotent Slodowy slices, including the subregular case
in type $A_2$.
Their framework fixes the graded Poisson central fibre.
The present work addresses the admissibility of
$R$-deformations of nonlinear Lie algebras and nonlinear
Lie conformal algebras, together with their PBW enveloping
constructions.
Comparing a finite quantization constructed here with
their universal family requires identifying its graded
Poisson central fibre and the corresponding parameter map.
Our chiral constructions further provide vertex algebra
quantizations of Poisson vertex algebras on the
corresponding arc spaces.

The quantum RDW spectral problem provides an
irregular oper realization related to the $A_2$
Toda chain, giving further evidence for the
conjectural connection between RDW-type integrable
systems and wild Hitchin system proposed
in~\cite{Fang26}. The geometric identification
of the corresponding integrable systems remains
an open problem. More generally, our constructions suggest further
questions concerning the quantization of nonlinear
$R$-deformations associated with other nilpotent
orbits and Slodowy slices.
In forthcoming work, we will develop a cohomological approach to
$R$-matrices and, in a separate paper~\cite{Fangpre}, use graded
symplectic geometry to study chiral $r$-matrices on Lie conformal
algebras, building on our joint work with N.~Ikeda~\cite{FI26}.

\section*{Acknowledgements}
The author would like to thank Professor Tomoyuki Arakawa
for suggesting the problem, for drawing his attention to the work
of De Sole and Kac after the author had completed his master's
thesis~\cite{Fang25}, and for reading the present manuscript.
 He is grateful to Professor Takahiro Shiota for valuable weekly discussions and for carefully reading the manuscript. He also thanks Professor Anton Alekseev for stimulating discussions and for his interest in the generalized AKS scheme during their meeting in Geneva in 2024. 
He also thanks Professor Anne Moreau for taking the time to listen to his explanation of the RDW system in Luminy. The author also acknowledges the use of OpenAI's GPT models for mathematical discussions and assistance in preparing the manuscript.

\subsection*{Grading conventions.}
Throughout this paper, unless otherwise stated, we work with the additive semigroup $\GP:=\frac12\mathbb Z_{>0}.$
All graded vector spaces and $\mathbb C[\partial]$-modules are understood to be $\Gamma_+$-graded, with $\partial$ and $\lambda$ assigned degree zero. The grading extends additively to tensor algebras, with the tensor unit having degree zero. Our constructions also extend to more general positively graded additive semigroups under suitable discreteness and homogeneity assumptions. We will also use notation $\GPz:=\GP\cup\{0\}=\frac12\ZZ_{\geq0}.$

\section{Preliminaries: Nonlinear Lie algebra and Nonlinear Lie conformal algebra}

Let $\fg=\bigoplus_{j\in\Gamma_+}\fg[j]$ be a $\Gamma_+$-graded
complex vector space equipped with a skew-symmetric bilinear map
\[
\mu:\fg\otimes\fg\longrightarrow T(\fg).
\]
For $a,b\in\fg$, write $[a,b]:=\mu(a,b)$.
For each nonzero homogeneous element $a\in\fg[j]$, set $\deg(a)=j$.
The grading on $\fg$ extends to a $\Gamma_{+,0}$-grading on $T(\fg)$
by setting
\[
\deg(1)=0,\qquad
\deg(a_1\otimes\cdots\otimes a_n)=\sum_{i=1}^n\deg(a_i)
\]
for $n\geq 1$ and nonzero homogeneous elements $a_1,\ldots,a_n\in\fg$.
Thus
\[
T(\fg)=\bigoplus_{r\in\Gamma_{+,0}}T(\fg)[r],
\qquad T(\fg)[0]=\CC 1.
\]
The induced filtration is given by
\[
T_{\leq j}(\fg):=
\bigoplus_{\substack{r\in\Gamma_{+,0}\\ r\leq j}}T(\fg)[r],
\qquad j\in\Gamma_{+,0}.
\]

We extend the bracket to bilinear maps
$\fg\otimes T(\fg)\longrightarrow T(\fg)$ and
$T(\fg)\otimes\fg\longrightarrow T(\fg)$ by the Leibniz rules
\[
\begin{aligned}
[a,A\otimes B] &= [a,A]\otimes B+A\otimes[a,B],\\
[A\otimes B,a] &= A\otimes[B,a]+[A,a]\otimes B,
\end{aligned}
\]
together with $[a,1]=[1,a]=0$, for all $a\in\fg$ and $A,B\in T(\fg)$.

For each $j\in\Gamma_{+,0}$, let $\mathcal M_j(\fg)$ be the subspace
of $T(\fg)$ spanned by the elements
\[
A\otimes\bigl(b\otimes c-c\otimes b-[b,c]\bigr)\otimes D,
\]
where $b\in\fg[p]$, $c\in\fg[q]$, $A\in T(\fg)[r]$, and
$D\in T(\fg)[s]$, with $p,q\in\Gamma_+$, $r,s\in\Gamma_{+,0}$,
and $r+p+q+s\leq j$. Set
\[
\mathcal M(\fg):=\bigcup_{j\in\Gamma_{+,0}}\mathcal M_j(\fg).
\]
Then $\mathcal M(\fg)$ is the two-sided ideal of $T(\fg)$ generated by
the elements $b\otimes c-c\otimes b-[b,c]$ for $b,c\in\fg$.
\begin{dfn}
The pair $(\fg,\mu)$ is called a \emph{nonlinear Lie algebra}
if the following conditions hold:
\begin{enumerate}
\item \textup{(Grading condition)}
For all $p,q\in\GP$, $\mu\bigl(\fg[p]\otimes\fg[q]\bigr)
\subseteq T_{\leq p+q-\frac12}(\fg).$
\item \textup{(Jacobi identity)}
For all $p,q,r\in\GP$ and all
$a\in\fg[p]$, $b\in\fg[q]$, and $c\in\fg[r]$,
\[
[a,[b,c]]-[b,[a,c]]-[[a,b],c]
\in \cM_{p+q+r-\frac12}(\fg).
\]
\end{enumerate}
\end{dfn}

\begin{thm}\label{thm:PBWNLA}
Let $(\fg,\mu)$ be a $\GP$-graded nonlinear Lie algebra.
Its universal enveloping algebra is the quotient
$U(\fg):=T(\fg)/\cM(\fg)$.
Let $\pi\colon T(\fg)\to U(\fg)$ be the quotient map.
The induced associative product satisfies
$\pi(A)\pi(B)=\pi(A\otimes B)$ for all $A,B\in T(\fg)$,
and its unit is $1_{U(\fg)}=\pi(1)$.

Let $\{x_i\}_{i\in\mathcal I}$ be a homogeneous basis of $\fg$
indexed by a totally ordered set $\mathcal I$.
Then $1_{U(\fg)}$ and the monomials
$\pi(x_{i_1})\cdots\pi(x_{i_n})$ with $n\geq 1$ and
$i_1\leq\cdots\leq i_n$ form a $\mathbb C$-basis of $U(\fg)$.
\end{thm}
\begin{exa}\label{exa:subW3}
Let $V=\CC u_1\oplus\CC u_2\oplus\CC u_3\oplus\CC\chi$
be a purely even $\GP$-graded vector space with normalized
Kazhdan degrees
\[
\deg u_1=\deg u_3=\frac32,\qquad
\deg u_2=1,\qquad \deg\chi=2.
\]
Define $\mu\colon V\otimes V\to T(V)$ by skew-symmetry and the relations
\[
\begin{alignedat}{2}
[u_1,u_2] &= 2u_1, &\qquad [u_2,u_3] &= 2u_3,\\
[u_3,u_1] &= u_2\otimes u_2+\chi, & [\chi,x] &= 0 \quad (x\in V).
\end{alignedat}
\]
For all $p,q\in\GP$, these relations give
\[
\mu\bigl(V[p]\otimes V[q]\bigr)
\subseteq T(V)[p+q-1]
\subseteq T_{\leq p+q-\frac12}(V).
\]

The Jacobi identity holds identically in $T(V)$ for all $a,b,c\in V$.
Indeed, by skew-symmetry and the centrality of $\chi$, it suffices
to check the triple $(u_1,u_2,u_3)$. Put $Q:=u_2\otimes u_2+\chi$.
Since $[u_2,Q]=0$, we obtain
\[
[u_1,[u_2,u_3]]-[u_2,[u_1,u_3]]-[[u_1,u_2],u_3]
=-2Q+0+2Q=0.
\]
Thus $(V,\mu)$ is a nonlinear Lie algebra. With $H=u_2$, $E=u_3$, $F=u_1$, and $C=\chi$, the defining
relations of $U(V,\mu)$ give the standard presentation of the
quantum finite $W$-algebra $\cW^\fin(\fsl_3,f_\sub)$,
where $f_\sub=E_{21}$. Hence
$U(V,\mu)\cong\cW^\fin(\fsl_3,f_\sub)$.
\end{exa}

We now recall the definition of a nonlinear Lie conformal algebra.

Let $V=\bigoplus_{j\in\GP}V[j]$ be a $\GP$-graded
$\CC[\partial]$-module with $\partial V[j]\subseteq V[j]$
for all $j\in\GP$. Write $\deg(a)=j$ for $0\ne a\in V[j]$.
As above, equip the tensor algebra
$T(V)=\bigoplus_{r\geq0}V^{\otimes r}$ over $\CC$ with the total
grading, with $T(V)[0]=V^{\otimes0}=\CC1$, and write
$T_{\leq j}(V)=\bigoplus_{k\leq j}T(V)[k]$ for its associated
filtration, where $j,k\in\GPz$.

\begin{dfn}\label{dfn:NCA}
A \emph{nonlinear conformal algebra} is a $\CC[\partial]$-module
endowed with a $\GP$-grading by $\CC[\partial]$-submodules
\[
V=\bigoplus_{p\in\GP}V[p],
\]
together with a $\CC$-linear map $[\,\cdot_\lambda\cdot\,]:
V\otimes V\longrightarrow T(V)[\lambda],
\ 
T(V)[\lambda]:=\CC[\lambda]\otimes T(V),$
satisfying the conditions below.
Here $\partial$ is extended to $T(V)$ as a derivation:
\[
\partial1=0,\qquad
\partial(A\otimes B)
=(\partial A)\otimes B+A\otimes(\partial B),
\qquad A,B\in T(V),
\]
and acts coefficientwise on $T(V)[\lambda]$.

\begin{enumerate}
\item \textup{(Sesquilinearity)}
For all $a,b\in V$, $[\partial a_\lambda b]=-\lambda[a_\lambda b],
\
[a_\lambda\partial b]=(\lambda+\partial)[a_\lambda b].$

\item \textup{(Grading condition)}
For all $p,q\in\GP$, $[V[p]_\lambda V[q]]
\subseteq T_{\leq p+q-\frac12}(V)[\lambda].$
\end{enumerate}
\end{dfn}
Let $N:T(V)\otimes T(V)\longrightarrow T(V),
\ 
L_\lambda:T(V)\otimes T(V)\longrightarrow T(V)[\lambda]$
be the maps defined by the recursive construction
of~\cite[Lemma~3.2]{DSK05}. In particular,
\[
N(a,B)=a\otimes B,\qquad
L_\lambda(a,b)=[a_\lambda b],
\qquad a,b\in V,\quad B\in T(V).
\]

For each $j\in\GPz$, let $\cM_j(V)$ be the $\CC$-linear span
of all elements
\[
A\otimes\left(
(b\otimes c-c\otimes b)\otimes D
-
N\left(
\int_{-\partial}^{0}L_\lambda(b,c)\,d\lambda,D
\right)
\right),
\]
where $b\in V[p],\  c\in V[q],\ 
A\in T(V)[r],\ D\in T(V)[s],$
with $p,q\in\GP$, $r,s\in\GPz$, and $r+p+q+s\leq j$.
Set
\[
\cM(V):=\bigcup_{j\in\GPz}\cM_j(V).
\]
The formal integral is defined coefficientwise by
\[
\int_{-\partial}^{0}
\left(\sum_{k\geq0}\lambda^k F_k\right)d\lambda
:=
\sum_{k\geq0}\frac{(-1)^k}{k+1}\partial^{k+1}F_k,
\qquad F_k\in T(V),
\]
where only finitely many $F_k$ are nonzero.
Thus, in the definition of $\cM_j(V)$, these derivatives
act on the coefficients of $L_\lambda(b,c)$ before
applying $N(\,\cdot\,,D)$.

\begin{dfn}\label{dfn:NLCA}
A \emph{nonlinear Lie conformal algebra} is a nonlinear conformal
algebra $V$ whose $\lambda$-bracket satisfies the following
additional axioms.
\begin{enumerate}
\item \textup{(Skewsymmetry)}
For all $a,b\in V$, we have $[a_\lambda b]=-[b_{-\lambda-\partial}a]$,
where $\partial$ acts on the coefficients of the bracket.

\item \textup{(Filtered Jacobi identity)}
For all $p,q,r\in\GP$ and all $a\in V[p]$, $b\in V[q]$, $c\in V[r]$,
\[
L_\lambda(a,L_\mu(b,c))
-L_\mu(b,L_\lambda(a,c))
-L_{\lambda+\mu}(L_\lambda(a,b),c)
\in \CC[\lambda,\mu]\otimes\cM_{p+q+r-\frac12}(V).
\]
\end{enumerate}
\end{dfn}
\subsection*{Notation}
For $a,b\in V$, we also write $\mu_\lambda(a,b)$ for $[a_\lambda b]$.
\section{filtered $R$-matrices of nonlinear Lie algebra and splitting constructions}
In this section, we introduce filtered $R$-matrices for nonlinear Lie algebras and establish a splitting theorem for such $R$-matrices. We also provide several examples illustrating this construction.

\begin{dfn}\label{dfn:FRM}
Let $(\fg,\mu)$ be a $\GP$-graded nonlinear Lie algebra, and let
$R\in\End_\CC(\fg)$ satisfy $R(\fg[j])\subseteq\fg[j]$ for every
$j\in\GP$. For $a,b\in\fg$, set
$[a,b]_R:=\mu(Ra,b)+\mu(a,Rb)$, and denote the resulting
skew-symmetric map $\fg\otimes\fg\to T(\fg)$ by $\mu_R$.
We call $R$ a \emph{filtered $R$-matrix} on $(\fg,\mu)$ if
$(\fg,\mu_R)$ is a nonlinear Lie algebra with respect to the
original $\GP$-grading.
\end{dfn}
\begin{dfn}\label{dfn:SC}
Let $(V,\mu)$ be a $\GP$-graded nonlinear Lie algebra, and let
$R:V\longrightarrow V$ be a filtered $R$-matrix of $(V,\mu)$.
Recall that $\mu_R(a,b):=\mu(R(a),b)+\mu(a,R(b)),
\ a,b\in V.$
We say that $R$ is \emph{strongly compatible} with $\mu$ if,
for every $t\in\CC$, the bracket $\mu_t:=\mu+t\mu_R:V\otimes V\longrightarrow T(V)$
defines a nonlinear Lie algebra structure on $V$ with respect
to the original $\GP$-grading.
\end{dfn}
\begin{dfn}\label{dfn:STRJ}
Let $(V,\mu)$ be a $\GP$-graded nonlinear Lie algebra.
Define its Jacobiator $\Jac(\mu):V^{\otimes3}\longrightarrow T(V)$
by $\Jac(\mu)(a,b,c)
:=[a,[b,c]]-[b,[a,c]]-[[a,b],c],
\ a,b,c\in V.$
We say that $\mu$ satisfies the \emph{strict Jacobi identity}
if $\Jac(\mu)(a,b,c)=0
\ \text{in }T(V)$
for all $a,b,c\in V.$

The Jacobi axiom in the definition of a nonlinear Lie algebra
will be referred to as the \emph{graded Jacobi identity}.
Explicitly, it requires $\Jac(\mu)(a,b,c)
\in\cM_{p+q+r-\frac12}^{\mu}(V)$
for all $p,q,r\in\GP$ and all
$a\in V[p]$, $b\in V[q]$, and $c\in V[r]$,
where $\cM_j^{\mu}(V)$ denotes the relation space associated
with $\mu$.
\end{dfn}

\begin{dfn}\label{dfn:SPT}
Let $(\fg,\mu)$ be a $\GP$-graded nonlinear Lie algebra.
A \emph{splitting} of $(\fg,\mu)$ is a decomposition $\fg=\fg_+\oplus\fg_-$
into graded subspaces such that $\mu(\fg_\pm\otimes\fg_\pm)\subseteq T(\fg_\pm).$
The associated \emph{splitting operator} is the
degree preserving $\CC$-linear map $R:=P_+-P_-,$
where $P_\pm:\fg\longrightarrow\fg$ are the projections
defined by $P_\pm(a_++a_-):=a_\pm,
\ a_\pm\in\fg_\pm.$
If $R$ is a filtered $R$-matrix of $(\fg,\mu)$, we call it
a \emph{splitting $R$-matrix}.
\end{dfn}

\begin{lemm}\label{lemm:SUBSP}
Let $(\fg,\mu)$ be a $\GP$-graded nonlinear Lie algebra,
and let $\fh\subseteq\fg$ be a graded subspace such that $\mu(\fh\otimes\fh)\subseteq T(\fh).$
Denote the restricted bracket by $\mu_{\fh}:\fh\otimes\fh\longrightarrow T(\fh).$
Then, for every $j\in\GPz$, $\cM_j^\mu(\fg)\cap T(\fh)
=
\cM_j^{\mu_{\fh}}(\fh).$
Consequently, $(\fh,\mu_{\fh})$ is a nonlinear Lie algebra
with respect to the induced grading $\fh=\bigoplus_{p\in\GP}\fh[p],
\ 
\fh[p]:=\fh\cap\fg[p].$
\end{lemm}

\begin{proof}
Fix $j\in\GPz$. By definition,
$\cM_j^{\mu_{\fh}}(\fh)\subseteq\cM_j^\mu(\fg)\cap T(\fh)$.
For the reverse inclusion, take $F\in\cM_j^\mu(\fg)\cap T(\fh)$
and extend a homogeneous ordered basis of $\fh$ to one of $\fg$.

Since $F\in T_{\leq j}(\fh)$, we may reorder its monomials inside
$T(\fh)$ using $xy\equiv yx+\mu_{\fh}(x,y)$ for $x>y$.
Each swap reduces the inversion number, and each bracket term
has lower weight. Thus $F\equiv Q\pmod{\cM_j^{\mu_{\fh}}(\fh)}$
for some linear combination $Q$ of ordered monomials in $\fh$.
Since $F$ maps to zero in $U(\fg,\mu)$, so does $Q$.
These ordered monomials belong to the PBW basis of $U(\fg,\mu)$
after projection, so $Q=0$. This proves the equality.

The restricted bracket inherits skewsymmetry and the grading condition.
For $a\in\fh[p]$, $b\in\fh[q]$, and $c\in\fh[r]$, its Jacobiator
agrees with that of $\mu$ and lies in
$\cM_{p+q+r-\frac12}^\mu(\fg)\cap T(\fh)$.
The equality just proved gives the Jacobi condition for $\mu_{\fh}$.
\end{proof}

\begin{thm}\label{thm:SPR}
Let $\fg=\fg_+\oplus\fg_-$ be a splitting as in
Definition~\ref{dfn:SPT}, and put
\[
\fg_\pm[p]:=\fg_\pm\cap\fg[p],
\qquad p\in\GP.
\]
Suppose that there exists $\delta\in\GP$ such that
\[
\mu(\fg_\pm[p]\otimes\fg_\pm[q])
\subseteq T(\fg_\pm)[p+q-\delta]
\]
for all $p,q\in\GP$, where $T(\fg_\pm)[m]$ denotes
the homogeneous component of total weight $m$,
with the convention that $T(\fg_\pm)[m]=0$ for $m<0$.
Then the associated splitting operator $R=P_+-P_-$
is a filtered $R$-matrix of $(\fg,\mu)$.
In particular, the bracket
\[
\mu_R(a,b):=\mu(R(a),b)+\mu(a,R(b)),
\qquad a,b\in\fg,
\]
defines a nonlinear Lie algebra structure on $\fg$
with respect to the original $\GP$-grading.
\end{thm}

\begin{proof}
Since the splitting is graded, $R=P_+-P_-$ preserves the grading.
Hence $\mu_R$ is skew-symmetric and satisfies the grading condition.
Writing $\mu_\pm$ for the restrictions of $\mu$ to
$\fg_\pm\otimes\fg_\pm$, we have
\[
\mu_R|_{\fg_+\otimes\fg_+}=2\mu_+,\qquad
\mu_R|_{\fg_-\otimes\fg_-}=-2\mu_-,\qquad
\mu_R|_{\fg_+\otimes\fg_-}=0.
\]
By Lemma~\ref{lemm:SUBSP}, both $(\fg_+,\mu_+)$ and
$(\fg_-,\mu_-)$ are nonlinear Lie algebras.

We first show that $(\fg_+,2\mu_+)$ and $(\fg_-,-2\mu_-)$
are nonlinear Lie algebras. Take $(H,\nu,\kappa)$ to be either
$(\fg_+,\mu_+,2)$ or $(\fg_-,\mu_-,-2)$.
Since $2\delta$ is a positive integer, we may choose
$\zeta\in\CC^\times$ such that $\zeta^{2\delta}=\kappa$.
Define $S:T(H)\to T(H)$ by setting
$S(F)=\zeta^{2m}F$ for $F\in T(H)[m]$ and extending linearly.
The additivity of the grading makes $S$ a degree preserving
algebra automorphism. For $x\in H[p]$ and $y\in H[q]$,
the homogeneity hypothesis gives
\[
S\bigl(\kappa\nu(x,y)\bigr)
=\kappa\zeta^{2(p+q-\delta)}\nu(x,y)
=\zeta^{2(p+q)}\nu(x,y)
=\nu(Sx,Sy).
\]
Thus $S$ sends each defining relation for $\kappa\nu$ to the
corresponding relation for $\nu$ with the same degree bound.
Since $S$ is invertible and preserves each homogeneous component,
\[
S\bigl(\cM_j^{\kappa\nu}(H)\bigr)=\cM_j^\nu(H),
\qquad j\in\GPz.
\]

Now take $a\in H[p]$, $b\in H[q]$, and $c\in H[r]$, and set
$j=p+q+r-\frac12$. By the Leibniz rules, $S$ also intertwines
the iterated brackets, so
\[
S\bigl(\Jac(\kappa\nu)(a,b,c)\bigr)
=\Jac(\nu)(Sa,Sb,Sc)\in\cM_j^\nu(H).
\]
Applying $S^{-1}$ gives
$\Jac(\kappa\nu)(a,b,c)\in\cM_j^{\kappa\nu}(H)$.
Thus $(H,\kappa\nu)$ is a nonlinear Lie algebra.

It remains to verify the Jacobi condition for $\mu_R$.
By trilinearity, it suffices to consider homogeneous triples
$a,b,c$ with each entry in one of the two summands.
If all three lie in $\fg_\pm$, the preceding argument gives
\[
\Jac(\mu_R)(a,b,c)
\in\cM_j^{\pm2\mu_\pm}(\fg_\pm)
\subseteq\cM_j^{\mu_R}(\fg),
\qquad j=\deg(a)+\deg(b)+\deg(c)-\frac12.
\]
For a mixed triple, the vanishing of $\mu_R(\fg_+,\fg_-)$
and the Leibniz rules give
\[
[T(\fg_\pm),\fg_\mp]_R=0,\qquad
[\fg_\mp,T(\fg_\pm)]_R=0.
\]
Hence every term in the Jacobiator vanishes.
This proves that $R$ is a filtered $R$-matrix.
\end{proof}
\begin{rmk}\label{rmk:SPR}
As in the classical case, a splitting $R$-matrix $R$ on a
nonlinear Lie algebra $(V,\mu)$ satisfies $R^2=\Id_V$.
The corresponding splitting is recovered from $R$ as $V=V_+\oplus V_-,
\
V_\pm=\ker(R\mp\Id_V),$
with projections $P_\pm=\frac12(\Id_V\pm R).$
Moreover, since $R$ is assumed to preserve the grading,
both eigenspaces are graded subspaces:
\[
V_\pm=\bigoplus_p V_\pm[p],
\qquad
V_\pm[p]=V_\pm\cap V[p].
\]
We refer to this graded eigenspace decomposition as the
\emph{splitting property} of $R$.
It provides useful constraints for the classification of
splitting $R$-matrices.
\end{rmk}
\begin{thm}\label{thm:subW3SP}
Let $(V,\mu)$ be the nonlinear Lie algebra in
Example~\ref{exa:subW3}.
With rows recording the images of the ordered basis
$(u_1,u_2,u_3,\chi)$, the degree preserving splitting-type
$R$-matrices are precisely
\[
R=\pm\Id_V,\qquad
R=\pm R_1(t,\epsilon),\qquad
R=\pm R_2(t,\epsilon),
\]
where
\[
R_1(t,\epsilon)=
\begin{bmatrix}
1&0&0&0\\
0&1&0&0\\
-2t&0&-1&0\\
0&0&0&\epsilon
\end{bmatrix},
\qquad
R_2(t,\epsilon)=
\begin{bmatrix}
-1&0&-2t&0\\
0&1&0&0\\
0&0&1&0\\
0&0&0&\epsilon
\end{bmatrix},
\]
with $t\in\CC$ and $\epsilon\in\{1,-1\}$.
The operators $\pm\Id_V$ correspond to the trivial
splittings. Thus, up to replacing $R$ by $-R$, there
are exactly two families of nontrivial degree preserving
splitting type $R$-matrices.
\end{thm}

\begin{proof}
Let $R$ be a degree preserving involution of $V$.
Since $V[1]=\CC u_2$, we may replace $R$ by $-R$ and assume
that $R(u_2)=u_2$. Using the brackets in Example~\ref{exa:subW3},
one checks that the eigenspaces $V_\pm=\ker(R\mp\Id_V)$ satisfy
$\mu(V_\pm\otimes V_\pm)\subseteq T(V_\pm)$ precisely when
\[
R=\Id_V,\qquad R=R_1(t,\epsilon),\qquad\text{or}\qquad
R=R_2(t,\epsilon),
\qquad t\in\CC,\quad \epsilon\in\{\pm1\}.
\]
For each such splitting, the restricted brackets inherit the
homogeneity condition with $\delta=1$. Hence
Theorem~\ref{thm:SPR} applies to both $R$ and $-R$, giving
the asserted list.
\end{proof}
More generally, one can verify that the following operators define
mixed-type $R$-matrices for the subregular finite $\mathcal{W}$-algebras
associated with $\mathfrak{sl}_n$ and the Levy's algebra of type $D_4$.
The corresponding RDW-type integrable systems and their chiralizations
will be investigated in forthcoming work~\cite{FangGeo}. 

The following two propositions follow from direct computations.
\begin{prop}\label{prop:subWsln}
Let $n\geq3$, and choose Smith generators
\[
E,\ H,\ F,\ c_2,\ldots,c_{n-1}
\]
of the finite subregular $W$-algebra
$U(\mathfrak{sl}_n,f_{\mathrm{sub}})$ as in
\cite[Section~4]{GK}, where $c_2,\ldots,c_{n-1}$ are central.
Let
\[
V=\CC E\oplus\CC H\oplus\CC F
  \oplus\bigoplus_{j=2}^{n-1}\CC c_j
\]
carry the corresponding nonlinear Lie algebra structure,
with the normalized Kazhdan grading $\deg H=1,\ 
\deg E=\deg F=\frac n2,\
\deg c_j=j.$

Consider degree preserving linear maps $R:V\to V$
of mixed type, namely,
\[
R(E),R(F)\in\CC E\oplus\CC F,\qquad
R(H)\in\CC H,\qquad
R(c_j)\in\CC c_j.
\]
With rows recording the images of the ordered basis
$(E,H,F,c_2,\ldots,c_{n-1})$, all such filtered
$R$-matrices are
\[
R_\pm(a,b,p,r;\boldsymbol{\rho})
=
\begin{pmatrix}
a&0&p&0\\
0&b&0&0\\
r&0&\pm a&0\\
0&0&0&D_{\boldsymbol{\rho}}
\end{pmatrix},
\qquad
D_{\boldsymbol{\rho}}
=\operatorname{diag}(\rho_2,\ldots,\rho_{n-1}),
\]
where $a,b,p,r,\rho_2,\ldots,\rho_{n-1}\in\CC.$
The parameters $\rho_j$ do not affect the deformed bracket.
Moreover, the strongly compatible mixed-type
$R$-matrices are precisely those in the family $R_+$.
\end{prop}

We use Levy's presentation of the $D_4$-type algebra
$H(P,\gamma)$, with
$P(u)=3u^2+\alpha u+\beta$
\cite[Definition~1.1 and Section~3]{Le09}.

\begin{prop}\label{prop:subD4mix}
Consider the symmetric nonlinear Lie presentation of Levy's $D_4$-type algebra
\cite{Le09}.
Let
\[
Z=\CC\alpha\oplus\CC\beta\oplus\CC\gamma,
\qquad
V=\CC u\oplus\CC v\oplus\CC w\oplus Z,
\]
with the grading
\[
\deg u=\deg v=\deg\alpha=2,\qquad
\deg w=3,\qquad
\deg\beta=\deg\gamma=4.
\]
The nonlinear bracket $\mu:V\otimes V\to T(V)$ is
determined by skew-symmetry and
\begin{align*}
\mu(u,v)&=2w,\ 
\mu(u,w)=-(uv+vu)+\gamma,\\
\mu(v,w)&=v^2+3u^2+\alpha u+\beta,\ 
\mu(Z,V)=0,
\end{align*}
where tensor product symbols are omitted.

Consider degree preserving linear maps $R:V\to V$
of mixed type, namely,
\[
R(u),R(v)\in\CC u\oplus\CC v,\qquad
R(w)\in\CC w,\qquad R(Z)\subseteq Z.
\]
With rows recording the images of the ordered basis
$(u,v,w,\alpha,\beta,\gamma)$, all such filtered
$R$-matrices are
\[
R_{\mathrm{I}}=
\begin{pmatrix}
a&b&0&0\\
c&-a&0&0\\
0&0&e&0\\
0&0&0&D
\end{pmatrix},
\qquad
R_{\mathrm{II}}=
\begin{pmatrix}
a&0&0&0\\
0&a&0&0\\
0&0&e&0\\
0&0&0&D
\end{pmatrix},
\]
where $a,b,c,e\in\CC$, and $D=\rho\oplus M,\ 
\rho\in\CC,\  M\in\operatorname{Mat}_2(\CC),$
is an arbitrary degree preserving endomorphism of $Z$.
The block $D$ does not affect the deformed bracket.
Moreover, the strongly compatible mixed-type
$R$-matrices are precisely those in the family
$R_{\mathrm{II}}$.
\end{prop}

\begin{proof}
The assertions follow from a direct computation of the
generator Jacobiators for $\mu_R$ and $\mu+t\mu_R$,
together with the PBW theorem for nonlinear Lie algebras.
\end{proof}
\begin{rmk}
With $\alpha,\beta,\gamma$ retained as independent
central generators, the undeformed enveloping algebra
$U(V,\mu)$ in Proposition~\ref{prop:subD4mix} can be shown to be
isomorphic, as a filtered algebra, to the subregular
finite $W$-algebra
$U(\mathfrak{so}_8,f_{\mathrm{sub}})$,
see~\cite{Le09,ACET24}.
\end{rmk}
\subsection{Examples of filtered $R$-matrices with applications in RDW system}\label{Sec:EFMRDW}
In this section, we present several examples of filtered $R$-matrices. In particular, we apply our construction to the $R$-deformed classical finite $\mathcal W$-algebras introduced in our previous work~\cite{Fang26}, obtaining their PBW quantizations. Let $R$ be a classical $R$-matrix for the classical finite $\cW$-algebra in the sense of~\cite{Fang26}. We follow the conventions of Example~\ref{exa:subW3} and consider the vertical $R$-matrices studied in~\cite[Proposition~5.7]{Fang26}, where the action on the central element $\chi$ is defined by $R(\chi)=0$.  Following~\cite[Proposition~5.7]{Fang26}, we recall three families of vertical classical $R$-matrices for the subregular classical finite $\cW$-algebra of $\mathfrak{sl}_3$. We use the conventions of Example~\ref{exa:subW3}. Let
\[
V=\operatorname{Span}_{\CC}\{u_1,u_2,u_3,\chi\}.
\]
With respect to the ordered basis $(u_1,u_2,u_3,\chi)$, we represent $R\in\End_{\CC}(V)$ by
\[
R=
\begin{pmatrix}
R_{11}&R_{12}&R_{13}&0\\
R_{21}&R_{22}&R_{23}&0\\
R_{31}&R_{32}&R_{33}&0\\
0&0&0&0
\end{pmatrix},
\]
where $R(u_i)=\sum_{j=1}^3R_{ij}u_j$ and $R(\chi)=0$.

\begin{prop}\label{prop:RDWR}
Let $A:=\mathcal O(\Sf)=\CC[u_1,u_2,u_3,\chi]$ be the Poisson algebra of the subregular Slodowy slice of $\mathfrak{sl}_3$, with the Poisson bracket as in Example~\ref{exa:subW3}. Let $R\in\End_{\CC}(V)$ satisfy $R(\chi)=0$ and $R(\operatorname{Span}_{\CC}\{u_1,u_2,u_3\})\subseteq\operatorname{Span}_{\CC}\{u_1,u_2,u_3\}$. Define
\[
\{u_i,u_j\}_R:=\{R(u_i),u_j\}+\{u_i,R(u_j)\},
\qquad \{\chi,-\}_R=0,
\]
and extend the bracket to $A$ by the Leibniz rules. Then the following three families give vertical classical $R$-matrices on $\Sf$.
\begin{description}
\item[Type I (Cartan type)]
Assume $R$ is diagonal. Then $R$ defines a Poisson bracket on $A$ if and only if
\[
R=\operatorname{diag}(a,b,\pm a,0),
\qquad a,b\in\CC.
\]

\item[Type II (upper triangular type)]
Assume $R$ is upper triangular, i.e., $R_{21}=R_{31}=R_{32}=0$. Then the Jacobi identity is equivalent to
\[
R_{11}^2=R_{33}^2,\qquad
R_{12}(R_{11}+R_{22})=0,\qquad
R_{23}(R_{22}+R_{33})=0.
\]

\item[Type III (mixed type)]
Assume $R_{12}=R_{21}=R_{23}=R_{32}=0$. Then $R$ defines a Poisson bracket on $A$ if and only if
\[
R_{11}=\pm R_{33}.
\]
\end{description}
In particular, every Type~III$_+$ $R$-matrix, characterized by $R_{11}=R_{33}$, defines a compatible Poisson pencil
\[
\{-,-\}_t:=\{-,-\}+t\{-,-\}_R,
\qquad t\in\CC,
\]
on the Slodowy slice $\Sf$.
\end{prop}
We work with the full subregular Slodowy slice in $\mathfrak{sl}_3$,
retaining the central generator $\chi$. The parameter $\rho\in\CC$ is arbitrary and does not affect
the deformed bracket.

The following theorem gives a PBW quantization of the
$R$-deformed Poisson algebra on $\Sf$.
Taking the quotient by the central generator $\chi$
also gives a PBW quantization of its central reduction
at $\chi=0$.
\begin{thm}\label{thm:PBWQRDW}
Let
\[
V=\CC u_1\oplus\CC u_2\oplus\CC u_3\oplus\CC\chi
\]
be the vector space from Example~\ref{exa:subW3}, equipped
with the $\GP$-grading
\[
V[1]=\CC u_2,\qquad
V\!\left[\frac32\right]=\CC u_1\oplus\CC u_3,\qquad
V[2]=\CC\chi,
\]
with all other homogeneous components equal to zero.
Define a skew-symmetric $\CC$-linear map
\[
\mu:V\otimes V\longrightarrow T(V)
\]
by
\[
\begin{aligned}
\mu(u_1,u_2)&=2u_1,\ 
\mu(u_2,u_3)=2u_3,\\
\mu(u_1,u_3)&=-(u_2\otimes u_2+\chi),\ 
\mu(\chi,x)=0,\  x\in V.
\end{aligned}
\]
Then $(V,\mu)$ is a nonlinear Lie algebra.

For any $a,b,d,p,r,\rho\in\CC$ satisfying $a^2=d^2$,
let $R\in\operatorname{End}_{\CC}(V)$ be given by
\[
R=
\begin{pmatrix}
a&0&p&0\\
0&b&0&0\\
r&0&d&0\\
0&0&0&\rho
\end{pmatrix},
\]
where rows record the images of the ordered basis
$(u_1,u_2,u_3,\chi)$.
Then $R$ is a filtered $R$-matrix of $(V,\mu)$.
In particular, the deformed bracket
\[
\mu_R(x,y):=\mu(R(x),y)+\mu(x,R(y)),
\qquad x,y\in V,
\]
defines a nonlinear Lie algebra structure on $V$
with respect to the same $\GP$-grading.

Moreover, $R$ is strongly compatible with $\mu$ if and
only if $a=d$. The universal enveloping algebra
$U_R:=U(V,\mu_R)$ is a filtered PBW quantization of
the Slodowy slice $\Sf$ equipped with the
$R$-deformed Poisson bracket.
\end{thm}

\begin{proof}
By Example~\ref{exa:subW3}, $(V,\mu)$ is a nonlinear Lie algebra
and $\mu$ has degree $-1$. Since $R$ preserves each homogeneous
component of $V$, the bracket $\mu_R$ also has degree $-1$ and
hence satisfies the grading condition. Set $Q:=u_2\otimes u_2+\chi$.
A direct computation gives
\[
\begin{aligned}
\mu_R(u_1,u_2)&=2(a+b)u_1-2pu_3, &
\mu_R(u_1,u_3)&=-(a+d)Q,\\
\mu_R(u_2,u_3)&=2(b+d)u_3-2ru_1, &
\mu_R(\chi,x)&=0 \qquad (x\in V).
\end{aligned}
\]

For $s,t\in\CC$, put $\nu_{s,t}:=s\mu+t\mu_R$.
By skew-symmetry and the centrality of $\chi$, it suffices to
compute the Jacobiator on $(u_1,u_2,u_3)$. Using $a^2=d^2$,
we obtain
\[
\Jac(\nu_{s,t})(u_1,u_2,u_3)
=2t(a-d)\bigl(s+t(a+d)\bigr)Q
=2st(a-d)Q.
\]
Taking $(s,t)=(0,1)$ shows that $\mu_R$ satisfies the Jacobi
identity in $T(V)$, so $R$ is a filtered $R$-matrix.
If $a=d$, the same identity holds for $\mu+t\mu_R$ for every
$t\in\CC$, proving strong compatibility. Conversely, if $R$ is
strongly compatible with $\mu$, then
$\Jac(\mu+\mu_R)(u_1,u_2,u_3)=2(a-d)Q$ has zero image in
$U(V,\mu+\mu_R)$. By the PBW theorem, $u_2^2$ and $\chi$ are
distinct basis elements in this algebra, so the image of $Q$ is
nonzero. Hence $a=d$.

Finally, let $U_\hbar$ be the $\CC[\hbar]$-algebra defined by
the relations $xy-yx=\hbar\mu_R(x,y)$ for $x,y\in V$.
The ordered monomials span $U_\hbar$ by reordering and are
linearly independent by PBW after specialization at every
$\hbar\in\CC$. Thus $U_\hbar$ is free over $\CC[\hbar]$, with
\[
U_\hbar/\hbar U_\hbar\cong\CC[u_1,u_2,u_3,\chi],
\qquad
U_\hbar/(\hbar-1)U_\hbar\cong U_R.
\]
The bracket induced by $\hbar^{-1}[-,-]$ on the first fiber
is the $R$-deformed Poisson bracket. Hence $U_R=U(V,\mu_R)$
is the asserted filtered PBW quantization. Since $\chi$ is
central, its reduction at $\chi=c$ is
$\widetilde{U}(V,\mu_R):=U_R/(\chi-c)$ for each $c\in\CC$.
\end{proof}

\section{filtered R-matrix of nonlinear Lie conformal algebra}
In this section, we introduce and study filtered $R$-matrices for nonlinear Lie conformal algebras, which we also call \emph{chiralized filtered $R$-matrices}. We establish a chiral analogue of Theorem~\ref{thm:SPR}.
\begin{dfn}\label{dfn:CFRM}
Let $(V,[\,\cdot_\lambda\cdot\,],\partial)$ be a
$\GP$-graded nonlinear Lie conformal algebra, and let \(R\)
be a degree preserving $\CC[\partial]$-linear map, that is, $R\partial=\partial R,
\
R(V[j])\subseteq V[j],
\ j\in\GP.$
Define the $\CC$-linear map
\[
[\,\cdot_\lambda\cdot\,]_R:
V\otimes_{\CC}V
\longrightarrow
\CC[\lambda]\otimes_{\CC}T(V)
\]
by
\[
[a_\lambda b]_R
:=[R(a)_\lambda b]+[a_\lambda R(b)],
\qquad a,b\in V.
\]
We say that $R$ is a \emph{chiral filtered $R$-matrix}
of $(V,[\,\cdot_\lambda\cdot\,],\partial)$ if
$(V,[\,\cdot_\lambda\cdot\,]_R,\partial)$ is a nonlinear
Lie conformal algebra with respect to the original
$\GP$-grading.
\end{dfn}
\begin{dfn}\label{dfn:CSC}
Let $(V,[\,\cdot_\lambda\cdot\,],\partial)$ be a
$\GP$-graded nonlinear Lie conformal algebra, and let
$R$ be a chiral filtered $R$-matrix on $V$.
Recall that
\[
[a_\lambda b]_R
:=[R(a)_\lambda b]+[a_\lambda R(b)],
\qquad a,b\in V.
\]
We say that $R$ is \emph{strongly compatible}
with $[\,\cdot_\lambda\cdot\,]$ if, for every $t\in\CC$,
the $\CC$-linear map
\[
[\,\cdot_\lambda\cdot\,]_t:
V\otimes_{\CC}V
\longrightarrow \CC[\lambda]\otimes_{\CC}T(V),
\]
defined by
\[
[a_\lambda b]_t
:=[a_\lambda b]+t[a_\lambda b]_R,
\qquad a,b\in V,
\]
makes $(V,[\,\cdot_\lambda\cdot\,]_t,\partial)$
a nonlinear Lie conformal algebra with respect to
the original $\GP$-grading.
Here the recursively defined maps
$N^{(t)}$, $L_\lambda^{(t)}$, and the relation spaces
$\cM_j^{(t)}(V)$ are those associated with
$[\,\cdot_\lambda\cdot\,]_t$.
\end{dfn}
\begin{dfn}\label{dfn:CSTRJ}
Let $(V,[\,\cdot_\lambda\cdot\,],\partial)$ be a
$\GP$-graded nonlinear conformal algebra whose
$\lambda$-bracket satisfies conformal skewsymmetry.
Let
\[
L_\lambda:
T(V)\otimes_{\CC}T(V)
\longrightarrow \CC[\lambda]\otimes_{\CC}T(V)
\]
be the canonical recursively defined extension of
$[\,\cdot_\lambda\cdot\,]$, and let
$\cM_j(V)$, $j\in\GPz$, denote the associated
degree bounded relation spaces.

For $a,b,c\in V$, set
\[
\begin{aligned}
\Jac_{\lambda,\nu}(a,b,c)
:=L_\lambda(a,L_\nu(b,c))-L_\nu(b,L_\lambda(a,c))-L_{\lambda+\nu}(L_\lambda(a,b),c),
\end{aligned}
\]
where the nested expressions are interpreted by
linearity in the formal parameters.

We say that $[\,\cdot_\lambda\cdot\,]$ satisfies the
\emph{strict Jacobi identity} if
\[
\Jac_{\lambda,\nu}(a,b,c)=0
\quad\text{in }\CC[\lambda,\nu]\otimes_{\CC}T(V)
\]
for all $a,b,c\in V$.

We say that $[\,\cdot_\lambda\cdot\,]$ satisfies the
\emph{graded Jacobi identity} if
\[
\Jac_{\lambda,\nu}(a,b,c)
\in
\CC[\lambda,\nu]\otimes_{\CC}
\cM_{p+q+r-\frac12}(V)
\]
for all $p,q,r\in\GP$ and all
$a\in V[p]$, $b\in V[q]$, and $c\in V[r]$.
\end{dfn}

\begin{dfn}\label{dfn:CSPT}
Let $(V,[\,\cdot_\lambda\cdot\,],\partial)$ be a
$\GP$-graded nonlinear Lie conformal algebra over $\CC$.
A \emph{splitting} of
$(V,[\,\cdot_\lambda\cdot\,],\partial)$
is a decomposition $V=V_+\oplus V_-$
into graded $\CC[\partial]$-submodules such that $[V_\pm{}_\lambda V_\pm]
\subseteq
\CC[\lambda]\otimes_{\CC}T(V_\pm).$
Here $T(V_\pm)$ is identified with its natural image
in $T(V)$.

The associated \emph{splitting operator} is the
degree preserving $\CC[\partial]$-linear map $R=P_+-P_-$
where $P_\pm:V\longrightarrow V$ are the projections defined by $P_\pm(a_++a_-)=a_\pm,
\ a_\pm\in V_\pm.$
In particular, $R\partial=\partial R,
\ 
R(V[j])\subseteq V[j],
\ j\in\GP.$
If, in addition, $R$ is a chiral filtered $R$-matrix
on $V$, we call it a \emph{chiral splitting $R$-matrix}.
\end{dfn}

We now give a general construction of chiral splitting
$R$-matrices. We begin with the chiral analogue of
Lemma~\ref{lemm:SUBSP}.

\begin{lemm}\label{lemm:CSUBSP}
Let $(V,\partial,[\,_\lambda\,])$ be a $\frac12\ZZ_{>0}$-graded nonlinear Lie conformal algebra, and let $H\subseteq V$ be a graded $\CC[\partial]$-submodule such that $[H_\lambda H]\subseteq\CC[\lambda]\otimes T(H).$
Denote the restricted $\lambda$-bracket by
\[
[\,_\lambda\,]_H:H\otimes H\longrightarrow\CC[\lambda]\otimes T(H).
\]
Then $(H,\partial,[\,_\lambda\,]_H)$ is a nonlinear Lie conformal algebra with respect to the induced grading.
\end{lemm}
\begin{proof}
Fix $j\in\GPz$. The recursive definitions of $N$ and $L_\lambda$
restrict to $T(H)$, so $\cM_j(H)\subseteq\cM_j(V)\cap T(H)$.
For the reverse inclusion, take $F\in\cM_j(V)\cap T(H)$ and
extend a homogeneous ordered $\CC$-basis of $H$ to one of $V$.

Since $F\in T_{\leq j}(H)$, we may reorder its monomials inside
$T(H)$. For adjacent basis elements $x>y$, the defining relation
gives, omitting tensor product symbols,
\[
AxyD\equiv AyxD+
A\,N_H\left(\int_{-\partial}^{0}[x_\lambda y]_H\,d\lambda,D\right)
\pmod{\cM_j(H)}
\]
whenever $AxyD$ has weight at most $j$. Each swap reduces the
inversion number, and the correction term has lower weight.
Thus $F\equiv Q\pmod{\cM_j(H)}$ for some linear combination $Q$
of ordered monomials in $H$. Since $F$ maps to zero in
$T(V)/\cM(V)$, so does $Q$. The ordered monomials in $Q$ project
to distinct PBW basis elements, so $Q=0$. This proves the equality.

The restricted bracket inherits sesquilinearity, skewsymmetry,
and the grading condition. For $a\in H[p]$, $b\in H[q]$, and
$c\in H[r]$, its Jacobiator agrees with that in $V$. Each
coefficient therefore lies in
$\cM_{p+q+r-\frac12}(V)\cap T(H)
=\cM_{p+q+r-\frac12}(H)$, as required.
\end{proof}
\begin{thm}\label{thm:CSPT}
Let $(V,\mu_\lambda,\partial)$ be a $\GP$-graded
nonlinear Lie conformal algebra, and let
\[
V=V_+\oplus V_-
\]
be a splitting as in Definition~\ref{dfn:CSPT}.
Put
\[
V_\pm[p]:=V_\pm\cap V[p],
\qquad p\in\GP.
\]
Suppose that there exists a fixed $\delta\in\GP$ such that
\[
\mu_\lambda(V_\pm[p],V_\pm[q])
\subseteq
\CC[\lambda]\otimes_{\CC}T(V_\pm)[p+q-\delta]
\]
for all $p,q\in\GP$, where $T(V_\pm)[m]$ denotes
the homogeneous component of total weight $m$,
with $T(V_\pm)[m]=0$ for $m<0$.

Then the associated splitting operator
\[
R=P_+-P_-
\]
is a chiral filtered $R$-matrix of
$(V,\mu_\lambda,\partial)$.
In particular, the bracket
\[
\mu_{R,\lambda}(a,b)
:=\mu_\lambda(R(a),b)+\mu_\lambda(a,R(b)),
\qquad a,b\in V,
\]
defines a nonlinear Lie conformal algebra structure on
$V$ with respect to the original $\GP$-grading.
\end{thm}

\begin{proof}
Since the splitting is graded and $V_\pm$ are
$\CC[\partial]$-submodules, $R=P_+-P_-$ preserves the grading
and commutes with $\partial$. Hence $\mu_{R,\lambda}$ satisfies
sesquilinearity, conformal skewsymmetry, and the grading condition.
Its restrictions to $V_+\otimes V_+$ and $V_-\otimes V_-$ are
$2\mu_{+,\lambda}$ and $-2\mu_{-,\lambda}$, respectively,
and its mixed brackets vanish.

Take $(H,\nu_\lambda,\kappa)$ to be either
$(V_+,\mu_{+,\lambda},2)$ or $(V_-,\mu_{-,\lambda},-2)$.
By Lemma~\ref{lemm:CSUBSP}, $(H,\nu_\lambda,\partial)$ is a
nonlinear Lie conformal algebra. As in the proof of
Theorem~\ref{thm:SPR}, choose $\zeta\in\CC^\times$ such that
$\zeta^{2\delta}=\kappa$ and define $S:T(H)\to T(H)$ by
$S(F)=\zeta^{2m}F$ for $F\in T(H)[m]$. By homogeneity,
\[
S\bigl(\kappa\nu_\lambda(a,b)\bigr)=\nu_\lambda(Sa,Sb),
\qquad a,b\in H,
\]
where $S$ acts coefficientwise. Since $S$ preserves tensor
products and commutes with $\partial$, the recursive formulas
show that $S$ also intertwines $N$ and $L_\lambda$ for the two
brackets. Thus $S$ identifies their relation spaces at each
degree bound, and the Jacobi condition for $\kappa\nu_\lambda$
follows from that for $\nu_\lambda$.

For triples in one summand, the Jacobi condition now follows
from $\cM_j^{\pm2\mu_\pm}(V_\pm)\subseteq\cM_j^{\mu_R}(V)$.
For mixed triples, the left and right Wick formulas give
\[
L^{\mu_R}_\lambda(V_\pm,T(V_\mp))
=L^{\mu_R}_\lambda(T(V_\pm),V_\mp)=0.
\]
Hence every mixed Jacobiator vanishes. By trilinearity, this
proves the Jacobi condition for $\mu_{R,\lambda}$.
\end{proof}
\section{PBW quantizations of $R$-deformed Poisson vertex algebras}\label{sec:PBWPVA}

In this section, we construct explicit nonlinear Lie conformal
algebras whose universal enveloping vertex algebras provide
PBW quantizations of $R$-deformed Poisson vertex algebras.
In particular, we extend the genuinely first-order
deformations in~\cite[Proposition~5.18]{Fang26} to the
arc algebra of the subregular Slodowy slice $\Sf$
in $\fsl_3$ and construct flat families of vertex algebras
quantizing them, retaining $\chi$ as an independent
central differential generator. Following the conventions of Section~\ref{Sec:EFMRDW}, let $\partial$
denote the derivation of the jet algebra and
consider the $\mathbb C[\partial]$-linear operators defined by
\[
\begin{aligned}
R_4(u_1)&=R_4(u_2)=R_4(\chi)=0,
& R_4(u_3)&=au_1+b\partial u_1,\qquad b\neq0,\\
R_7(u_2)&=R_7(u_3)=R_7(\chi)=0,
& R_7(u_1)&=Bu_3+D\partial u_3,\qquad D\neq0.
\end{aligned}
\]
Define their deformed brackets on generators by
\[
\{x_\lambda y\}_{R_i}
:=\{R_i(x)_\lambda y\}_0+\{x_\lambda R_i(y)\}_0,
\qquad
x,y\in\{u_1,u_2,u_3,\chi\},
\]
and extend by sesquilinearity and the PVA Leibniz rules.
For $R_4$, the nonzero independent brackets are
\[
\begin{aligned}
\{u_2{}_\lambda u_3\}_{R_4}
=-2\bigl(a+b(\lambda+\partial)\bigr)u_1,\ 
\{u_3{}_\lambda u_3\}_{R_4}
=b(\partial+2\lambda)(u_2^2+\chi).
\end{aligned}
\]
Here $u_1$ and $\chi$ are central.
For $R_7$, they are
\[
\begin{aligned}
\{u_1{}_\lambda u_2\}_{R_7}
=-2(B-D\lambda)u_3,\ 
\{u_1{}_\lambda u_1\}_{R_7}
=-D(\partial+2\lambda)(u_2^2+\chi),
\end{aligned}
\]
and $u_3$ and $\chi$ are central.
In both cases $\{u_2{}_\lambda u_2\}_{R_i}=0$.
The following proposition proves, in particular,
that these formulas define PVAs.

\begin{prop}\label{prop:firstorder-Slodowy-PBW}
The $R_4$- and $R_7$-deformed PVAs on $\mathcal A$
admit PBW quantizations over $\CC[\hbar]$. More precisely, put $A:=\CC[\hbar]$ and let
\[
R_\hbar
=A[\partial]u_1\oplus A[\partial]u_2
 \oplus A[\partial]u_3\oplus A[\partial]\chi.
\]
All generators are even. Equip $R_\hbar$
with the $\GP$-grading
\[
\deg u_1=\deg u_3=\frac32,\qquad
\deg u_2=1,\qquad \deg\chi=2,
\qquad
\deg\partial=\deg\lambda=\deg\hbar=0.
\]
Set $Q:=u_2\otimes u_2+\chi\in T_A(R_\hbar).$
For $i=4,7$, define an $A$-linear map
\[
[\,\cdot_\lambda\cdot\,]_{i,\hbar}:
R_\hbar\otimes_AR_\hbar
\longrightarrow
A[\lambda]\otimes_A T_A(R_\hbar)
\]
by the following formulas, together with sesquilinearity
and conformal skewsymmetry:
\[
\begin{gathered}
[u_2{}_\lambda u_3]_{4,\hbar}
=-2\hbar\bigl(a+b(\lambda+\partial)\bigr)u_1,\qquad
[u_3{}_\lambda u_3]_{4,\hbar}
=b\hbar(\partial+2\lambda)Q,\\
[u_1{}_\lambda x]_{4,\hbar}
=[\chi_\lambda x]_{4,\hbar}
=[u_2{}_\lambda u_2]_{4,\hbar}=0
\quad(x\in R_\hbar);
\end{gathered}
\]
and
\[
\begin{gathered}
[u_1{}_\lambda u_2]_{7,\hbar}
=-2\hbar(B-D\lambda)u_3,\qquad
[u_1{}_\lambda u_1]_{7,\hbar}
=-D\hbar(\partial+2\lambda)Q,\\
[u_3{}_\lambda x]_{7,\hbar}
=[\chi_\lambda x]_{7,\hbar}
=[u_2{}_\lambda u_2]_{7,\hbar}=0
\quad(x\in R_\hbar).
\end{gathered}
\]
These formulas define nonlinear Lie conformal algebras
$R_{i,\hbar}$ over $A$ with \([R_\hbar[p]_\lambda R_\hbar[q]]_{i,\hbar}
\subseteq
A[\lambda]\otimes_A T_A(R_\hbar)[p+q-1].\)
The nonlinear universal enveloping vertex algebras $U_{i,\hbar}:=U_A(R_{i,\hbar})$
are free $A$-modules. Moreover, \(U_{i,\hbar}/\hbar U_{i,\hbar}
\cong
\bigl(\mathcal A,\{\,\cdot_\lambda\cdot\,\}_{R_i}\bigr).\)

The independent singular OPEs
in the $R_4$-branch are
\[
\begin{aligned}
H(z)F(w)
&\sim
-\frac{2b\hbar Z(w)}{(z-w)^2}
-\frac{2\hbar(aZ+b\partial Z)(w)}{z-w},\\
F(z)F(w)
&\sim
\frac{2b\hbar(:HH:+\chi)(w)}{(z-w)^2}
+\frac{b\hbar\,\partial(:HH:+\chi)(w)}{z-w}.
\end{aligned}
\]
In the $R_7$-branch, writing $E=u_1$, $H=u_2$,
and $Z=u_3$, they are
\[
\begin{aligned}
E(z)H(w)
&\sim
\frac{2D\hbar Z(w)}{(z-w)^2}
-\frac{2B\hbar Z(w)}{z-w},\\
E(z)E(w)
&\sim-
\frac{2D\hbar(:HH:+\chi)(w)}{(z-w)^2}
-\frac{D\hbar\,\partial(:HH:+\chi)(w)}{z-w}.
\end{aligned}
\]
In both cases all OPEs involving $Z$ or $\chi$,
as well as $H(z)H(w)$, are regular.
The reversed OPEs are determined by conformal skewsymmetry.

\end{prop}

\begin{proof}
We prove the assertion for $R_4$ and write $Z=u_1,\  H=u_2,\  F=u_3.$
Sesquilinearity, skewsymmetry, and the grading condition
follow directly from the defining brackets.

Since $Z$ and $\chi$ are central and
$[H_\lambda H]_\hbar=0$, all Jacobiators on
generators vanish except possibly the triple $(F,F,F)$.
By the Wick formulas, we obtain
\[
\begin{aligned}
\Jac_{\lambda,\eta}(F,F,F)
=2b\hbar^2(\lambda-\eta)
\bigl(a+b(\lambda+\eta)+b\partial\bigr)\cdot(H\otimes Z-Z\otimes H).
\end{aligned}
\]
Since $[H_\lambda Z]_\hbar=0$, we have $H\otimes Z-Z\otimes H
\in\cM_{5/2}^{(4,\hbar)}(R).$
As $\partial$ preserves the grading, it follows that
\[
\Jac_{\lambda,\eta}(F,F,F)
\in A[\lambda,\eta]\otimes_A
\cM_{5/2}^{(4,\hbar)}(R)
\subseteq A[\lambda,\eta]\otimes_A
\cM_4^{(4,\hbar)}(R),
\]
where $4=3\deg F-\frac12$.
Hence $\mathcal R_{4,\hbar}$ is a nonlinear Lie
conformal algebra.

The case of $R_7$ follows by replacing
$(Z,H,F)$ with $(u_3,u_2,u_1)$ and
$(a,b)$ with $(-B,-D)$.

Finally, by the PBW theorem~\cite[Theorem~3.9]{DSK05},
$U_{i,\hbar}$ is free over $A=\CC[\hbar]$, and
\[
U_{i,\hbar}/\hbar U_{i,\hbar}
\simeq
S\bigl(\CC[\partial]u_1\oplus\CC[\partial]u_2
\oplus\CC[\partial]u_3\oplus\CC[\partial]\chi\bigr).
\]
Moreover, the semiclassical $\lambda$-bracket
induced by $\hbar^{-1}[\cdot_\lambda\cdot]_\hbar$
coincides with the corresponding $R_i$-deformed
PVA bracket. This completes the proof.
\end{proof}

\begin{rmk}
Proposition~\ref{prop:firstorder-Slodowy-PBW} constructs
PBW quantizations of the classical $R_i$-deformed PVAs.
It does not by itself identify $R_i$ as a chiral
filtered $R$-matrix of a specified undeformed
nonlinear Lie conformal algebra; that assertion
requires a separate verification for the chosen
quantum bracket.
\end{rmk}

\section{Applications to integrable systems}\label{Sec:APPIST}
In this section, we apply the Generalized AKS scheme~\cite{Fang25} to obtain several examples of integrable systems in classical and quantum sense. Our main example is the RDW system introduced in~\cite{Fang26}.
\subsection{Classical RDW system}\label{sec:CRDW}
We first recall the \(R\)-deformed system associated with the subregular nilpotent element of \(\mathfrak{sl}_3\), introduced in~\cite{Fang26} on the Slodowy slice, and provide some geometric background. Further algebro-geometric aspects of this system will be studied in~\cite{FangGeo}.
We now describe the RDW system associated with
Theorem~\ref{thm:PBWQRDW}, extending the construction
of~\cite[Theorem~5.8]{Fang26} to the Slodowy slice
$\Sf$.

Adopt the notation of Theorem~\ref{thm:PBWQRDW}.
Set $d=a$ and $\delta:=a-b$.
Recall the Hamiltonian
\[
S:=u_1u_3+\frac16u_2^3+\frac12\chi u_2.
\]
The corresponding RDW system on the Slodowy slice
$\Sf$, defined by $\dot u_i=\{S,u_i\}_R$,
takes the form
\begin{equation}\label{eq:GRDW}
\begin{alignedat}{2}
\dot u_1 &= (u_2^2+\chi)(\delta u_1+p u_3), \qquad&
\dot u_2 &= 2(r u_1^2-p u_3^2),\\
\dot u_3 &= -(u_2^2+\chi)(r u_1+\delta u_3), &
\dot\chi &= 0.
\end{alignedat}
\end{equation}
The polynomial $J:=r u_1^2+p u_3^2+2\delta u_1u_3$
is a first integral, while $K:=4aS-J$
and $\chi$ are Casimirs of the $R$-deformed Poisson
bracket. The parameter $\rho$ does not affect
the Poisson structure or the resulting dynamics.

In particular, for $a=b=d$ and $p=r=1$,
the system reduces to
\begin{equation}\label{eq:RDW}
\begin{alignedat}{2}
\dot u_1 &= u_3(u_2^2+\chi), \qquad&
\dot u_2 &= 2(u_1^2-u_3^2),\\
\dot u_3 &= -u_1(u_2^2+\chi), &
\dot\chi &= 0.
\end{alignedat}
\end{equation}
Its conserved quantities include
\[
S=u_1u_3+\frac16u_2^3+\frac12\chi u_2,
\qquad
J=u_1^2+u_3^2,
\qquad
\chi.
\]
The central reduction $\chi=0$ recovers the
CRDW system of~\cite[Theorem~5.8]{Fang26}.

We briefly discuss the geometry of the CRDW system; a more detailed
study will appear in~\cite{FangGeo}. In particular, its spectral curves
are related to the generalized Kummer surfaces studied by
Bertin--Vanhaecke~\cite{BV94}.

Fix $s,j\in\CC$ with $j(j^2-4s^2)\neq0$, and consider the common
level curve
\[
M_{s,j}:=
\left\{
(u_1,u_2,u_3)\in\AAA^3
\ \middle|\
u_1u_3+\frac16u_2^3=s,\quad u_1^2+u_3^2=j
\right\}.
\]
The functions $x=u_2,\ y=3(u_1^2-u_3^2)$
define a degree two morphism from $M_{s,j}$ to the affine curve
\[
C_{s,j}^{\mathrm{aff}}:
\quad y^2=-x^6+12sx^3+9(j^2-4s^2).
\]
More precisely, this morphism identifies $C_{s,j}^{\mathrm{aff}}$
with the quotient of $M_{s,j}$ by the involution $(u_1,u_2,u_3)\longmapsto(-u_1,u_2,-u_3).$
Under the above assumption, $C_{s,j}^{\mathrm{aff}}$ is smooth,
and its smooth projective model $C_{s,j}$ has genus two.
It admits an automorphism of order three, $\sigma:C_{s,j}\longrightarrow C_{s,j},
\ 
(x,y)\longmapsto(\zeta_3x,y),$
where $\zeta_3$ is a primitive cube root of unity.

After rescaling $x$ and $y$, the affine equation takes the form
\[
Y^2=X^6+2\kappa X^3+1,
\qquad
\kappa^2=\frac{4s^2}{4s^2-j^2}\neq1,
\]
studied in~\cite{BV94}. Let $\sigma_J$ denote the automorphism
of $\operatorname{Jac}(C_{s,j})$ induced by $\sigma$. Then $\operatorname{Jac}(C_{s,j})/\langle\sigma_J\rangle$
is the corresponding generalized Kummer surface.

We note that the RDW system is also related to Seiberg--Witten theory,
a connection we discussed in part in our previous paper~\cite{Fang26}.
The order-three automorphism $\sigma$ is induced by a generator of
$\mu_3\subset\CC^\times$ acting on $\Sf$ via the Kazhdan action.
For fixed $j\ne0$, the substitution $Y=iy$ gives
$Y^2=(x^3-6s)^2-9j^2$. As $s$ varies, this identifies the curves
$C_{s,j}$ with the $u=0$ subfamily of Seiberg--Witten curves for
pure $\mathcal N=2$ $\mathrm{SU}(3)$ gauge theory~\cite[Eq.~(4.4)]{KLT95}.
In~\cite{FangGeo}, we will study these connections further and
examine the relation between the RDW system and the even master
system of Bertin and Vanhaecke~\cite{BV94}.

\subsection{Quantization of RDW system}\label{sec:QRDW}
We use the PBW theorem to quantize the Poisson algebra underlying
our RDW system~\eqref{eq:GRDW}. We restrict to the subfamily $d=a$
in Theorem~\ref{thm:PBWQRDW}. Let $U_\hbar$ denote the
$\hbar$-adic completion of the $\mathbb C[\hbar]$-algebra
generated by $\hu_1,\hu_2,\hu_3,\hchi$, subject to the relations
\begin{equation}\label{eq:QRDW}
\begin{alignedat}{2}
[\hu_1,\hu_2] &= 2\hbar\bigl((a+b)\hu_1-p\hu_3\bigr), \qquad&
[\hu_2,\hu_3] &= 2\hbar\bigl((a+b)\hu_3-r\hu_1\bigr),\\
[\hu_1,\hu_3] &= -2a\hbar(\hu_2^2+\hchi), &
[\hchi,\hu_i] &= 0 \quad (i=1,2,3),
\end{alignedat}
\end{equation}
where $[A,B]:=AB-BA$.

\begin{prop}\label{prop:QRDW}
The algebra $U_\hbar$ is a flat formal PBW quantization of
the Poisson algebra underlying the RDW system~\eqref{eq:GRDW}.
Moreover, the elements
\[
\hat H
=\frac12[\hu_1,\hu_3]_+
+\frac16\hu_2^3
+\frac12\hchi\hu_2
+\frac{(a+b)^2-pr}{3}\hbar^2\hu_2
\]
and 
\[
\hat J
=
r\hu_1^2+p\hu_3^2
+(a-b)[\hu_1,\hu_3]_+,
\]
where $[A,B]_+:=AB+BA$, satisfy $[\hat H,\hat J]=0.$
Their images modulo $\hbar$ are the corresponding classical
Hamiltonians $S$ and $J$. The system can be written as:
\[
\begin{alignedat}{2}
\frac{d\hu_i(t)}{dt} &= \frac{1}{\hbar}[\hat H,\hu_i(t)], \qquad&
\hu_i(0) &= \hu_i \quad (i=1,2,3),\\
\frac{d\hchi(t)}{dt} &= 0, &
\hchi(0) &= \hchi.
\end{alignedat}
\]
where $i=1,2,3.$
\end{prop}

\begin{proof}
The PBW argument applies over $\mathbb C[\hbar]$.
Consequently, $U_\hbar$ is topologically free over
$\mathbb C[[\hbar]]$, and
$$
U_\hbar/\hbar U_\hbar
\cong \mathbb C[u_1,u_2,u_3,\chi],
$$
with the semiclassical bracket induced by
$\hbar^{-1}[\cdot,\cdot]$ coinciding with the
$R$-deformed Poisson bracket.

A direct calculation using the defining relations shows that
$\hat K:=4a\hat H-\hat J$ is central in $U_\hbar$,
and hence $[\hat H,\hat J]=0$.
Finally, reducing $\hat H$ and $\hat J$ modulo $\hbar$
recovers the classical Hamiltonians $S$ and $J$.
\end{proof}

Motivated by the conjecture proposed in~\cite{Fang26},
we further investigate the connection between the
quantized RDW system and Hitchin systems.

\begin{thm}[Quantum RDW--Toda correspondence]
\label{thm:RDW-Toda-oper}
Consider the quantum RDW system with $a=b=0,\  p=r=1,$
and fix the central characters
\[
\hat\chi=c,\qquad
\hat J=\hat u_1^2+\hat u_3^2=j,
\qquad c\in\CC,\quad j\in\CC^\times.
\]
Let $\hat H$ be the corresponding quantum Hamiltonian.
Put
\[
\epsilon=-4i\hbar,\qquad
\mu=\frac94j^2,\qquad
D_z=\epsilon z\partial_z.
\]
Then the quantum spectral problem admits a
differential-operator realization
\[
6\mathsf H_z
=
D_z^3+
\left(3c+\frac{\epsilon^2}{8}\right)D_z
+z+\frac{\mu}{z}.
\]
In particular, the eigenvalue equation
$\mathsf H_z\psi=e\psi$ becomes
\begin{equation}\label{eq:RDW-quantum-spectral}
\left[
D_z^3+
\left(3c+\frac{\epsilon^2}{8}\right)D_z
+z+\frac{\mu}{z}-6e
\right]\psi=0.
\end{equation}

Choose $\Lambda\in\CC^\times$ with $\Lambda^6=\mu$,
and set $\zeta=\frac{z}{i\Lambda^3},
\
\hbar_{\mathrm T}=i\epsilon=4\hbar,$
together with
\[
E_2=3c+\frac{\epsilon^2}{8},
\qquad
E_3=-6e,
\qquad
t_{\mathbf E}(s)=s^3+E_2s+E_3.
\]
Under the change of variables
$\varphi(\zeta)=\psi(i\Lambda^3\zeta)$,
equation~\eqref{eq:RDW-quantum-spectral}
takes the Toda form
\begin{equation}\label{eq:RDW-Toda-form}
t_{\mathbf E}
(-i\hbar_{\mathrm T}\zeta\partial_\zeta)
\varphi
=
\Lambda^3
\left(i^3\zeta+i^{-3}\zeta^{-1}\right)\varphi.
\end{equation}
This coincides with the $N=3$ quantum spectral
equation of the periodic Toda chain
in~\cite[Eq.~(3.2)]{BRT25}.

Moreover, the spectral equation defines an
$\mathrm{SL}_3$-oper on $\mathbb G_m$ over
$\CC((\epsilon))$, extending meromorphically
to $\mathbb P^1$ with irregular singularities
at $0$ and $\infty$.
Its semiclassical spectral curve is
\[
p^3+3cp+z+\frac{\mu}{z}-6e=0,
\qquad
\lambda=p\,\frac{dz}{z}.
\]
\end{thm}

\begin{proof}
The centrally reduced algebra admits a representation on
$M=\CC[t,t^{-1}][[\hbar]]$ given by $\rho(\hat u_1+i\hat u_3)=t,\ 
\rho(\hat u_1-i\hat u_3)=jt^{-1},\ 
\rho(\hat u_2)=-2i\hbar t\partial_t.$
The operator $\rho(\hat H)$ preserves
$M^+=\CC[t^2,t^{-2}][[\hbar]]$.
Under $z=\frac{3}{2i}t^2$, its restriction to $M^+$
gives the asserted operator $\mathsf H_z$ and
equation~\eqref{eq:RDW-quantum-spectral}.
The substitution $z=i\Lambda^3\zeta$ then yields
the Toda form~\eqref{eq:RDW-Toda-form}.

Since $t_E$ has no quadratic term, the corresponding
companion connection is trace-free and, equipped with
the standard oper flag, defines an $\mathrm{SL}_3$-oper
on $\mathbb G_m$ over $\CC((\epsilon))$.
It extends meromorphically to $\mathbb P^1$
with irregular singularities at $0$ and $\infty$.
Finally, taking the semiclassical symbol of
\eqref{eq:RDW-quantum-spectral} gives the stated spectral curve.
\end{proof}

\begin{rmk}
The Toda spectral equation in
Theorem~\ref{thm:RDW-Toda-oper}
is associated with the wild Hitchin system
of type $A_2$ on $\mathbb P^1$;
see~\cite[Section~1.4]{BRT25}.
This provides further evidence for the
conjectural relationship between RDW type
integrable systems and wild Hitchin system
proposed in~\cite{Fang26}.
\end{rmk}

\begin{rmk}
One can check that, as in the Section~\ref{sec:CRDW} under the change of variables
\[
p=x,\qquad
z=-\frac12(x^3-6s+iy),\qquad
\mu=\frac94j^2,
\]
the affine curve $C_{s,j}^{\mathrm{aff}}$ is isomorphic to
\[
p^3+z+\frac{\mu}{z}-6s=0.
\]
This coincides with the semiclassical spectral curve
of the quantum RDW Hamiltonian in Theorem~\ref{thm:RDW-Toda-oper}
at $\chi=0$.
\end{rmk}

\subsection{Chiralization of the RDW system and its quantization}

In this section, we study the chiralization and
quantization of the RDW system.

We adopt the notation of Theorem~\ref{thm:PBWQRDW}
and restrict to the mixed family $d=a$.
Put $\sigma:=a+b$ and $\Delta:=\sigma^2-pr$.

Consider the jet algebra
$$
\cA=\mathcal O(J_\infty\Sf)
=\CC[T^n u_1,T^n u_2,T^n u_3,T^n\chi
\mid n\geq0].
$$
Here $T$ denotes the canonical derivation of $\cA$.
Below, we use the same symbol for the derivation on the tensor
algebra and the translation operator on its vertex algebra quotient;
under the identification of the classical limit with $\cA$,
the induced derivation agrees with $T$.
The corresponding $R$-deformed Poisson bracket
extends to a PVA structure on $\cA$ by
sesquilinearity and the Leibniz rules.
We now construct its chiral PBW quantization.

\begin{prop}\label{prop:CRDWPBW}
For every $a,b,p,r\in\CC$, the mixed PVA $\cA$ admits a PBW quantization.
\end{prop}

\begin{proof}\label{Prop:mixPBA}
Let $H=\CC[T]u_1\oplus\CC[T]u_2\oplus\CC[T]u_3\oplus\CC[T]\chi$, with grading
\[
\deg u_1=\deg u_3=\frac32,\qquad
\deg u_2=1,\qquad \deg\chi=2,\qquad \deg T=\deg\lambda=0.
\]
Define $[\cdot_\lambda\cdot]:H\otimes H\to\CC[\lambda]\otimes T(H)$ by
\begin{equation}\label{eq:RDWmixedquantum}
\begin{gathered}
\begin{alignedat}{2}
[u_1{}_\lambda u_2] &= 2(\sigma u_1-pu_3), \qquad&
[u_2{}_\lambda u_3] &= -2(ru_1-\sigma u_3),\\
[u_1{}_\lambda u_1] &= -4ap(2\lambda+T)u_2, &
[u_3{}_\lambda u_3] &= -4ar(2\lambda+T)u_2,\\
[u_2{}_\lambda u_2] &= -4\Delta\lambda\mathbf1, &
[\chi{}_\lambda u_i] &= [\chi{}_\lambda\chi]=0 \quad(i=1,2,3),
\end{alignedat}\\
[u_1{}_\lambda u_3]
= -2a(u_2\otimes u_2+\chi)-4a\sigma(2\lambda+T)u_2-8a\Delta\lambda^2\mathbf1.
\end{gathered}
\end{equation}
All remaining brackets are determined by skewsymmetry and
sesquilinearity.

The grading condition follows immediately from the defining brackets.

Moreover, a direct calculation using the Wick formulas gives
\[
\operatorname{Jac}_{\lambda,\mu}(x,y,z)
\in\CC[\lambda,\mu]\otimes
\mathcal M_{\deg x+\deg y+\deg z-2}(H)
\]
for all generators $x,y,z\in\{u_1,u_2,u_3,\chi\}$.
Hence $H$ is a NLCA.

By the PBW theorem~\cite[Theorem~3.9]{DSK05}, the corresponding
$\hbar$-adic enveloping vertex algebra $U_\hbar$, obtained by
replacing $(a,b,p,r)$ with $(\hbar a,\hbar b,\hbar p,\hbar r)$
in the above brackets, satisfies
$$
U_\hbar/\hbar U_\hbar\cong\cA
$$
as Poisson vertex algebras, with the bracket induced by
$\hbar^{-1}[\cdot{}_\lambda\cdot]_\hbar$.
\end{proof}
\begin{rmk}\label{rmk:mixPBA}
Let $\mathcal V(a,b,p,r)$ be the universal enveloping
vertex algebra defined by the brackets in~\eqref{eq:RDWmixedquantum} of
Proposition~\ref{prop:CRDWPBW}, before the $\hbar$-rescaling.
If $ap\Delta\neq0$, a direct change of generators,
together with the PBW property, gives $\mathcal V(a,b,p,r)
\simeq
\mathcal W^{-3}(\mathfrak{sl}_3,f_{\mathrm{sub}});$
see the critical-level presentation
in~\cite[Section~4]{GK}.
The quotient by the vertex ideal generated by $\chi$
is the universal enveloping vertex algebra of
the nonlinear Lie conformal algebra $R^2_{-1}$
of Bakalov--De Sole~\cite[Section~2.2]{BDS09}.
\end{rmk}
\subsection{Chiral quantization of RDW system}
Following the vertex-algebraic approach to quantum integrable
hierarchies in~\cite{Wang26,FF96,FF11}, we adopt the following definition.

\begin{dfn}\label{dfn:QIST}
Let $(\cV,\partial,\{\cdot{}_\lambda\cdot\})$ be a Poisson
vertex algebra, and let
\[
\cH_n\in\cV/\partial\cV,\qquad n\in\ZZ_{\geq0},
\]
be the pairwise commuting Hamiltonians of a given classical
integrable hierarchy.

Let $V_\hbar$ be an $\hbar$-adic vertex algebra quantization
of $\cV$, with translation operator $T$. Assume that
\[
\mathcal F_\hbar:=V_\hbar/TV_\hbar
\]
is topologically free over $\CC[[\hbar]]$, and denote the
natural reduction map by
\[
\sigma:\mathcal F_\hbar\longrightarrow\cV/\partial\cV.
\]
Equip $\mathcal F_\hbar$ with the Lie bracket
\[
\left[\int a,\int b\right]_\hbar
:=\int b_{(0)}a,
\qquad a,b\in V_\hbar,
\]
where $\int a$ denotes the class of $a$ modulo $TV_\hbar$.

A \emph{quantization of the given integrable hierarchy}
is a pair
\[
\left(V_\hbar,
\{\widehat{\cH}_n\}_{n\in\ZZ_{\geq0}}\right),
\qquad
\widehat{\cH}_n\in\mathcal F_\hbar,
\]
such that
\[
\sigma(\widehat{\cH}_n)=\cH_n,
\qquad
[\widehat{\cH}_m,\widehat{\cH}_n]_\hbar=0
\quad
\text{for all }m,n\in\ZZ_{\geq0}.
\]
\end{dfn}
Choose a $T$-equivariant PBW isomorphism $\phi:\mathcal A[[\hbar]]\xrightarrow{\sim}V_\hbar$
of topological $\mathbb C[[\hbar]]$-modules. For instance, one may
take the symmetrization map defined on monomials in the jet variables by
\[
\phi(a_1\cdots a_n)
=\frac{1}{n!}\sum_{\sigma\in S_n}
N_n(a_{\sigma(1)},\ldots,a_{\sigma(n)}),
\]
where $N_n$ denotes the right-nested normally ordered product
and $\phi(1)=\mathbf{1}$.
It induces an isomorphism
\[
(\mathcal A/T\mathcal A)[[\hbar]]
\xrightarrow{\sim}V_\hbar/TV_\hbar.
\]
In particular, $V_\hbar/TV_\hbar$ is topologically free.

\begin{cor}\label{cor:CQRDW}
Let $S$ and $J$ be as above, and let
\[
\wS:=:S:,\qquad \wJ:=:J:
\]
denote their normally ordered PBW lifts to $V_\hbar$. Then
\[
\left[\int\wS,\int\wJ\right]=0
\qquad\text{in }V_\hbar/TV_\hbar.
\]
Thus $\int\wS$ and $\int\wJ$ form a commuting pair
of quantum Hamiltonians in the sense of Wang~\cite{Wang26},
lifting the classical local Hamiltonians
$\int S$ and $\int J$, respectively.
\end{cor}
\begin{proof}
Set $\widehat K:=4a\wS-\wJ$.
A direct computation using the Wick formulas,
quasicommutativity and quasiassociativity gives
\[
\widehat K_{(0)}u_i=0 \quad (i=1,2,3),
\qquad
\widehat K_{(0)}\chi=0.
\]
Since zero modes are continuous derivations commuting with $T$,
and $V_\hbar$ is topologically strongly generated by
$u_1,u_2,u_3,\chi$, it follows that
$\widehat K_{(0)}V_\hbar=0$.
Therefore,
\[
\begin{aligned}
\left[\int\wS,\int\wJ\right]_\hbar
=\int\wJ_{(0)}\wS=4a\int\wS_{(0)}\wS
-\int\widehat K_{(0)}\wS=0,
\end{aligned}
\]
where $\int\wS_{(0)}\wS=0$ by skewsymmetry modulo $TV_\hbar$.
\end{proof}
\begin{rmk}\label{rmk:CQRDW}
Naive normal ordering does not, in general, preserve the
involutivity of higher Hamiltonians. For example, specialize to $a=b=\frac12,\ p=r=1,$
and write $x=u_1$, $y=u_2$, and $z=u_3$.
Let $S,J$ and their normally ordered lifts $\wS,\wJ$
be as above. Put $A=\frac12(y^2+\chi)^2-8yxz,
\ 
\Omega_2
=
8A(S\,TJ-J\,TS)
+
16yJ(x^2-z^2)\,TS.$
A direct Wick calculation gives $\left[\int:\wS^2:,\int:\wJ^2:\right]
=
\hbar^2\int:\Omega_2:+O(\hbar^3),$
where $:\Omega_2:$ denotes a normally ordered PBW lift
of $\Omega_2$. Since Euler derivatives annihilate total derivatives, the identity
\[
\operatorname{coeff}_{T\chi}
\left(\frac{\delta\Omega_2}{\delta y}\right)
=
4J\bigl(
\chi^2-y^4+2y(x^2-z^2)+8yxz
\bigr)
\neq0
\]
shows that $\Omega_2\notin T\mathcal A$.
Under the $T$-equivariant PBW identification $V_\hbar/TV_\hbar
\simeq
(\mathcal A/T\mathcal A)[[\hbar]],$
the coefficient of $\hbar^2$ in the displayed commutator
is therefore nonzero.
Thus the classically commuting local Hamiltonians
$\int S^2$ and $\int J^2$ do not remain in involution
under these naive normally ordered lifts.
\end{rmk}


\end{document}